\documentstyle[11pt]{article}
\begin{document}
\newcommand{\p}{\parallel }
\makeatletter \makeatother
\newtheorem{th}{Theorem}[section]
\newtheorem{lem}{Lemma}[section]
\newtheorem{de}{Definition}[section]
\newtheorem{rem}{Remark}[section]
\newtheorem{cor}{Corollary}[section]
\renewcommand{\theequation}{\thesection.\arabic {equation}}

\title{{\bf Modular invariance and anomaly cancellation formulas in odd dimension}}

\author{ Kefeng Liu; Yong Wang\\
 }

\date{}
\maketitle

\begin{abstract}~~  By studying modular invariance
properties of some characteristic forms, we get some new
anomaly cancellation formulas on $(4r-1)$ dimensional manifolds. As an application, we derive some results on divisibilities of the index of Toeplitz operators on $(4r-1)$ dimensional spin
manifolds and some congruent formulas on characteristic number for $(4r-1)$ dimensional spin$^c$ manifolds. \\

\end{abstract}

\section{Introduction}
  \quad In 1983, the physicists Alvarez-Gaum\'{e} and Witten [AW]
  discovered the "miraculous cancellation" formula for gravitational
  anomaly which reveals a beautiful relation between the top
  components of the Hirzebruch $\widehat{L}$-form and
  $\widehat{A}$-form of a $12$-dimensional smooth Riemannian
  manifold. Kefeng Liu [Li1] established higher dimensional "miraculous cancellation"
  formulas for $(8k+4)$-dimensional Riemannian manifolds by
  developing modular invariance properties of characteristic forms.
  These formulas could be used to deduce some divisibility results. In
  [HZ1], [HZ2], for each $(8k+4)$-dimensional smooth Riemannian
  manifold, a more general cancellation formula that involves a
  complex line bundle was established. This formula was applied to
  ${\rm spin}^c$ manifolds, then an analytic Ochanine congruence
  formula was derived. In [CH], Qingtao Chen and Fei Han obtained more twisted cancellation formulas for $8k$ and $8k+4$ dimensional
  manifolds and they also applied their cancellation formulas to study divisibilities on spin manifolds and
  congruences on ${\rm spin}^c$ manifolds.\\
  \indent Another important application of modular invariance
properties of characteristic forms is to prove the rigidity theorem on elliptic genera. For example, see [Li, 2,3], [LM], [LMZ1,2]. For odd
   dimensional manifolds, we proved the similar rigidity theorem for elliptic genera under the condition that fixed point submanifolds are $1$-dimensional in [LW]. In [HY], Han and Yu dropped off our condition
   and proved more general odd dimensional rigidity theorem for elliptic genera. In [HY], in order to prove the rigidity theorem, they constructed some interesting modular forms under the condition
     that the $3$th de-Rham cohomology of manifolds
   vanishes. \\
   \indent In parallel, a natural question is whether we can get some interesting cancellation formulas in odd dimension by modular forms constructed in [HY]. In this paper, we will give the confirmative answer of this question.
That is,  by studying modular invariance
properties of some characteristic forms, we get some new
anomaly cancellation formulas on $(4r-1)$ dimensional manifolds. As an application, we derive some results on divisibilities on $(4r-1)$ dimensional spin
manifolds and congruences on $(4r-1)$ dimensional spin$^c$ manifolds. In [CH1], a cancellation formula on $11$-dimensional manifolds was derived.  To the authors' best knowledge, our cancellation formulas appear for the first time for general odd dimensional manifolds.\\
\indent This paper is organized as follows: In Section 2, we review
some knowledge on characteristic forms and modular forms that we are
going to use. In Section 3, we prove some odd dimensional cancellation
formulas and we apply them to get some results on divisibilities on the index of Toeplitz operators on
spin manifolds. In Section 4, we prove some odd cancellation
formulas involving a
  complex line bundle. By these formulas, we get some congruent formulas on characteristic number for odd spin$^c$ manifolds.\\

\section{Characteristic forms and modular forms}
 \quad The purpose of this section is to review the necessary knowledge on
characteristic forms and modular forms that we are going to use.\\

 \noindent {\bf  2.1 Characteristic forms. }Let $M$ be a Riemannian manifold.
 Let $\nabla^{ TM}$ be the associated Levi-Civita connection on $TM$
 and $R^{TM}=(\nabla^{TM})^2$ be the curvature of $\nabla^{ TM}$.
 Let $\widehat{A}(TM,\nabla^{ TM})$ and $\widehat{L}(TM,\nabla^{ TM})$
 be the Hirzebruch characteristic forms defined respectively by (cf.
 [Zh])
 $$\widehat{A}(TM,\nabla^{ TM})={\rm
 det}^{\frac{1}{2}}\left(\frac{\frac{\sqrt{-1}}{4\pi}R^{TM}}{{\rm
 sinh}(\frac{\sqrt{-1}}{4\pi}R^{TM})}\right),$$
 $$\widehat{L}(TM,\nabla^{ TM})={\rm
 det}^{\frac{1}{2}}\left(\frac{\frac{\sqrt{-1}}{2\pi}R^{TM}}{{\rm
 tanh}(\frac{\sqrt{-1}}{4\pi}R^{TM})}\right).\eqno(2.1)$$
   Let $F$, $F'$ be two Hermitian vector bundles over $M$ carrying
   Hermitian connection $\nabla^F,\nabla^{F'}$ respectively. Let
   $R^F=(\nabla^F)^2$ (resp. $R^{F'}=(\nabla^{F'})^2$) be the curvature of
   $\nabla^F$ (resp. $\nabla^{F'}$). If we set the formal difference
   $G=F-F'$, then $G$ carries an induced Hermitian connection
   $\nabla^G$ in an obvious sense. We define the associated Chern
   character form as
   $${\rm ch}(G,\nabla^G)={\rm tr}\left[{\rm
   exp}(\frac{\sqrt{-1}}{2\pi}R^F)\right]-{\rm tr}\left[{\rm
   exp}(\frac{\sqrt{-1}}{2\pi}R^{F'})\right].\eqno(2.2)$$
   For any complex number $t$, let
   $$\wedge_t(F)={\bf C}|_M+tF+t^2\wedge^2(F)+\cdots,~S_t(F)={\bf
   C}|_M+tF+t^2S^2(F)+\cdots$$
   denote respectively the total exterior and symmetric powers of
   $F$, which live in $K(M)[[t]].$ The following relations between
   these operations hold,
   $$S_t(F)=\frac{1}{\wedge_{-t}(F)},~\wedge_t(F-F')=\frac{\wedge_t(F)}{\wedge_t(F')}.\eqno(2.3)$$
   Moreover, if $\{\omega_i\},\{\omega_j'\}$ are formal Chern roots
   for Hermitian vector bundles $F,F'$ respectively, then
   $${\rm ch}(\wedge_t(F))=\prod_i(1+e^{\omega_i}t).\eqno(2.4)$$
   Then we have the following formulas for Chern character forms,
   $${\rm ch}(S_t(F))=\frac{1}{\prod_i(1-e^{\omega_i}t)},~
{\rm
ch}(\wedge_t(F-F'))=\frac{\prod_i(1+e^{\omega_i}t)}{\prod_j(1+e^{\omega_j'}t)}.\eqno(2.5)$$
\indent If $W$ is a real Euclidean vector bundle over $M$ carrying a
Euclidean connection $\nabla^W$, then its complexification $W_{\bf
C}=W\otimes {\bf C}$ is a complex vector bundle over $M$ carrying a
canonical induced Hermitian metric from that of $W$, as well as a
Hermitian connection $\nabla^{W_{\bf C}}$ induced from $\nabla^W$.
If $F$ is a vector bundle (complex or real) over $M$, set
$\widetilde{F}=F-{\rm dim}F$ in $K(M)$ or $KO(M)$.\\

\noindent{\bf 2.2 Some properties about the Jacobi theta functions
and modular forms}\\
   \indent We first recall the four Jacobi theta functions are
   defined as follows( cf. [Ch]):
   $$\theta(v,\tau)=2q^{\frac{1}{8}}{\rm sin}(\pi
   v)\prod_{j=1}^{\infty}[(1-q^j)(1-e^{2\pi\sqrt{-1}v}q^j)(1-e^{-2\pi\sqrt{-1}v}q^j)],\eqno(2.6)$$
$$\theta_1(v,\tau)=2q^{\frac{1}{8}}{\rm cos}(\pi
   v)\prod_{j=1}^{\infty}[(1-q^j)(1+e^{2\pi\sqrt{-1}v}q^j)(1+e^{-2\pi\sqrt{-1}v}q^j)],\eqno(2.7)$$
$$\theta_2(v,\tau)=\prod_{j=1}^{\infty}[(1-q^j)(1-e^{2\pi\sqrt{-1}v}q^{j-\frac{1}{2}})
(1-e^{-2\pi\sqrt{-1}v}q^{j-\frac{1}{2}})],\eqno(2.8)$$
$$\theta_3(v,\tau)=\prod_{j=1}^{\infty}[(1-q^j)(1+e^{2\pi\sqrt{-1}v}q^{j-\frac{1}{2}})
(1+e^{-2\pi\sqrt{-1}v}q^{j-\frac{1}{2}})],\eqno(2.9)$$ \noindent
where $q=e^{2\pi\sqrt{-1}\tau}$ with $\tau\in\textbf{H}$, the upper
half complex plane. Let
$$\theta'(0,\tau)=\frac{\partial\theta(v,\tau)}{\partial
v}|_{v=0}.\eqno(2.10)$$ \noindent Then the following Jacobi identity
(cf. [Ch]) holds,
$$\theta'(0,\tau)=\pi\theta_1(0,\tau)\theta_2(0,\tau)\theta_3(0,\tau).\eqno(2.11)$$
\noindent Denote $SL_2({\bf Z})=\left\{\left(\begin{array}{cc}
\ a & b  \\
 c  & d
\end{array}\right)\mid a,b,c,d \in {\bf Z},~ad-bc=1\right\}$ the
modular group. Let $S=\left(\begin{array}{cc}
\ 0 & -1  \\
 1  & 0
\end{array}\right),~T=\left(\begin{array}{cc}
\ 1 &  1 \\
 0  & 1
\end{array}\right)$ be the two generators of $SL_2(\bf{Z})$. They
act on $\textbf{H}$ by $S\tau=-\frac{1}{\tau},~T\tau=\tau+1$. One
has the following transformation laws of theta functions under the
actions of $S$ and $T$ (cf. [Ch]):
$$\theta(v,\tau+1)=e^{\frac{\pi\sqrt{-1}}{4}}\theta(v,\tau),~~\theta(v,-\frac{1}{\tau})
=\frac{1}{\sqrt{-1}}\left(\frac{\tau}{\sqrt{-1}}\right)^{\frac{1}{2}}e^{\pi\sqrt{-1}\tau
v^2}\theta(\tau v,\tau);\eqno(2.12)$$
$$\theta_1(v,\tau+1)=e^{\frac{\pi\sqrt{-1}}{4}}\theta_1(v,\tau),~~\theta_1(v,-\frac{1}{\tau})
=\left(\frac{\tau}{\sqrt{-1}}\right)^{\frac{1}{2}}e^{\pi\sqrt{-1}\tau
v^2}\theta_2(\tau v,\tau);\eqno(2.13)$$
$$\theta_2(v,\tau+1)=\theta_3(v,\tau),~~\theta_2(v,-\frac{1}{\tau})
=\left(\frac{\tau}{\sqrt{-1}}\right)^{\frac{1}{2}}e^{\pi\sqrt{-1}\tau
v^2}\theta_1(\tau v,\tau);\eqno(2.14)$$
$$\theta_3(v,\tau+1)=\theta_2(v,\tau),~~\theta_3(v,-\frac{1}{\tau})
=\left(\frac{\tau}{\sqrt{-1}}\right)^{\frac{1}{2}}e^{\pi\sqrt{-1}\tau
v^2}\theta_3(\tau v,\tau),\eqno(2.15)$$
 $$\theta'(v,\tau+1)=e^{\frac{\pi\sqrt{-1}}{4}}\theta'(v,\tau),~~
 \theta'(0,-\frac{1}{\tau})=\frac{1}{\sqrt{-1}}\left(\frac{\tau}{\sqrt{-1}}\right)^{\frac{1}{2}}
\tau\theta'(0,\tau).\eqno(2.16)$$
 \noindent {\bf Definition 2.1} A modular form over $\Gamma$, a
 subgroup of $SL_2({\bf Z})$, is a holomorphic function $f(\tau)$ on
 $\textbf{H}$ such that
 $$f(g\tau):=f\left(\frac{a\tau+b}{c\tau+d}\right)=\chi(g)(c\tau+d)^kf(\tau),
 ~~\forall g=\left(\begin{array}{cc}
\ a & b  \\
 c & d
\end{array}\right)\in\Gamma,\eqno(2.17)$$
\noindent where $\chi:\Gamma\rightarrow {\bf C}^{\star}$ is a
character of $\Gamma$. $k$ is called the weight of $f$.\\
Let $$\Gamma_0(2)=\left\{\left(\begin{array}{cc}
\ a & b  \\
 c  & d
\end{array}\right)\in SL_2({\bf Z})\mid c\equiv 0~({\rm
mod}~2)\right\},$$
$$\Gamma^0(2)=\left\{\left(\begin{array}{cc}
\ a & b  \\
 c  & d
\end{array}\right)\in SL_2({\bf Z})\mid b\equiv 0~({\rm
mod}~2)\right\},$$  be the two modular subgroups of $SL_2({\bf Z})$.
It is known that the generators of $\Gamma_0(2)$ are
$T,~ST^2ST$, the generators of $\Gamma^0(2)$ are $STS,~T^2STS$ (cf.[Ch]).\\
\indent If $\Gamma$ is a modular subgroup, let ${\mathcal{M}}_{{\bf
R}}(\Gamma)$ denote the ring of modular forms over $\Gamma$ with
real Fourier coefficients. Writing $\theta_j=\theta_j(0,\tau),~1\leq
j\leq 3,$ we introduce four explicit modular forms (cf. [Li1]),
$$\delta_1(\tau)=\frac{1}{8}(\theta_2^4+\theta_3^4),~~\varepsilon_1(\tau)=\frac{1}{16}\theta_2^4\theta_3^4,$$
$$\delta_2(\tau)=-\frac{1}{8}(\theta_1^4+\theta_3^4),~~\varepsilon_2(\tau)=\frac{1}{16}\theta_1^4\theta_3^4.\eqno(2.18)$$
\noindent They have the following Fourier expansions in
$q^{\frac{1}{2}}$:
$$\delta_1(\tau)=\frac{1}{4}+6q+6q^2\cdots,~~\varepsilon_1(\tau)=\frac{1}{16}-q+7q^2+\cdots,$$
$$\delta_2(\tau)=-\frac{1}{8}-3q^{\frac{1}{2}}-3q+\cdots,~~\varepsilon_2(\tau)=q^{\frac{1}{2}}+8q+\cdots,\eqno(2.19)$$
\noindent where the $"\cdots"$ terms are the higher degree terms,
all of which have integral coefficients. They also satisfy the
transformation laws,
$$\delta_2(-\frac{1}{\tau})=\tau^2\delta_1(\tau),~~~~~~\varepsilon_2(-\frac{1}{\tau})
=\tau^4\varepsilon_1(\tau).\eqno(2.20)$$

\noindent {\bf Lemma 2.2} ([Li1]) {\it $\delta_1(\tau)$ (resp.
$\varepsilon_1(\tau)$) is a modular form of weight $2$ (resp. $4$)
over $\Gamma_0(2)$, $\delta_2(\tau)$ (resp. $\varepsilon_2(\tau)$)
is a modular form of weight $2$ (resp. $4$) over $\Gamma^0(2)$ and
moreover ${\mathcal{M}}_{{\bf R}}(\Gamma^0(2))={\bf
R}[\delta_2(\tau),\varepsilon_2(\tau)]$.}

\section { Some cancellation formulas in odd dimension }

\quad Let $M$ be a $(4r-1)$ dimensional Riemannian manifold. Set
   $$\Theta_1(T_{C}M)=
   \bigotimes _{n=1}^{\infty}S_{q^n}(\widetilde{T_CM})\otimes
\bigotimes _{m=1}^{\infty}\wedge_{q^m}(\widetilde{T_CM})
,$$
$$\Theta_2(T_{C}M)=\bigotimes _{n=1}^{\infty}S_{q^n}(\widetilde{T_CM})\otimes
\bigotimes _{m=1}^{\infty}\wedge_{-q^{m-\frac{1}{2}}}(\widetilde{T_CM})$$
$$\Theta_3(T_{C}M)=\bigotimes _{n=1}^{\infty}S_{q^n}(\widetilde{T_CM})\otimes
\bigotimes _{m=1}^{\infty}\wedge_{q^{m-\frac{1}{2}}}(\widetilde{T_CM})
.\eqno(3.1)$$
\indent We recall the odd Chern character of a smooth map $g$ from $M$ to the general linear group $GL(N,{\bf{C}})$ with $N$ a positive integer (see [Zh]). Let $d$ denote a trivial connection on ${\bf{C}}^N|_M$.
We will denote by $c_n(M,[g])$ the cohomology class associated to the closed $n$-form
$$c_n({\bf C}^N|_M,g,d)=\left(\frac{1}{2\pi\sqrt{-1}}\right)^{\frac{(n+1)}{2}}{\rm Tr}[(g^{-1}dg)^n].\eqno(3.2)$$
The odd Chern character form ${\rm ch}({\bf C}^N|_M,g,d)$ associated to $g$ and $d$ by definition is
$${\rm ch}({\bf C}^N|_M,g,d)=\sum_{n=1}^{\infty}\frac{n!}{(2n+1)!}c_{2n+1}({\bf C}^N|_M,g,d).\eqno(3.3)$$
Let the connection $\nabla_u$ on the trivial bundle ${\bf C}^N|_M$ defined by
$$\nabla_u=(1-u)d+ug^{-1}\cdot d\cdot g,~~u\in [0,1].\eqno(3.4)$$
Then we have
$$d{\rm ch}({\bf C}^N|_M,g,d)={\rm ch}({\bf C}^N|_M,d)-{\rm ch}({\bf C}^N|_M,g^{-1}\cdot d\cdot g).\eqno(3.5)$$
\indent Now let $g:M\rightarrow SO(N)$ and we assume that $N$ is even and large enough. Let $E$ denote the trivial real vector bundle of rank $N$ over $M$. We equip $E$ with the canonical trivial metric and
trivial connection $d$. Set $$\nabla_u=d+ug^{-1}dg,~~u\in[0,1].$$ Let $R_u$ be the curvature of $\nabla_u$, then
$$R_u=(u^2-u)(g^{-1}dg)^2.\eqno(3.6)$$
We also consider the complexification of $E$ and $g$ extends to a unitary automorphism of $E_{\bf C}$. The connection $\nabla_u$ extends to a Hermitian connection on $E_{\bf C}$ with curvature still given by (3.6).
Let $\triangle(E)$ be the spinor bundle of $E$, which is a trivial Hermitian bundle of rank $2^{\frac {N}{2}}$. We assume that $g$ has a lift to the Spin group Spin$(N)$: $g^{\triangle}:M\rightarrow {\rm Spin}(N).$
So $g^{\triangle}$ can be viewed as an automorphism of $\triangle(E)$ preserving the Hermitian metric. We lift $d$ on $E$ to be a trivial Hermitian connection $d^\triangle$ on $\triangle(E)$, then
$$\nabla^\triangle_u=(1-u)d^\triangle+u(g^\triangle)^{-1}\cdot d^\triangle\cdot g^\triangle,~~u\in[0,1]\eqno(3.7)$$
lifts $\nabla^u$ on $E$ to $\triangle(E)$. Let $Q_j(E)$, $j=1,2,3$ be the virtual bundles defined as following:
 $$Q_1(E)=\triangle(E)\otimes
   \bigotimes _{n=1}^{\infty}\wedge_{q^n}(\widetilde{E_C})
,$$
$$Q_2(E)=\bigotimes _{n=1}^{\infty}\wedge_{-q^{n-\frac{1}{2}}}(\widetilde{E_C});~~
Q_3(E)=\bigotimes _{n=1}^{\infty}\wedge_{q^{n-\frac{1}{2}}}(\widetilde{E_C})
.\eqno(3.8)$$
Let $g$ on $E$ have a lift $g^{Q_j(E)}$ on $Q_j(E)$ and $\nabla_u$ have a lift $\nabla^{Q_j(E)}_u$ on $Q_j(E)$. Following [HY], we defined ${\rm ch}(Q_j(E),g^{Q_j(E)},d,\tau)$ for $j=1,2,3$ as following
$${\rm ch}(Q_j(E),\nabla^{Q_j(E)}_0,\tau)-{\rm ch}(Q_j(E),\nabla^{Q_j(E)}_1,\tau)=d{\rm ch}(Q_j(E),g^{Q_j(E)},d,\tau),\eqno(3.9)$$
where
$${\rm ch}(Q_1(E),g^{Q_1(E)},d,\tau)=-\frac{2^{N/2}}{8\pi^2}\int^1_0{\rm Tr}\left[g^{-1}dg\frac{\theta'_1(R_u/(4\pi^2),\tau)}{\theta_1(R_u/(4\pi^2),\tau)}\right]du,\eqno(3.10)$$
and for $j=2,3$
$${\rm ch}(Q_j(E),g^{Q_j(E)},d,\tau)=-\frac{1}{8\pi^2}\int^1_0{\rm Tr}\left[g^{-1}dg\frac{\theta'_j(R_u/(4\pi^2),\tau)}{\theta_j(R_u/(4\pi^2),\tau)}\right]du.\eqno(3.11)$$
By Proposition 2.2 in [HY], we have if $c_3(E_C,g,d)=0$, then for any integer $l\geq 1$ and $j=1,2,3$,~
${\rm ch}(Q_j(E),g^{Q_j(E)},d,\tau)^{(4l-1)}$ are modular forms of weight $2l$ over $\Gamma_0(2)$, $\Gamma^0(2)$ and $\Gamma_\theta$ respectively.
Let (see [HY, Def. 2.3])
$$\Phi_L(\nabla^{TM},g,d,\tau)=\widehat{L}(TM,\nabla^{TM}){\rm ch}(\Theta_1(TM),\nabla^{\Theta_1(TM)},\tau){\rm ch}(Q_1(E),g^{Q_1(E)},d,\tau);\eqno(3.12)$$
$$\Phi_W(\nabla^{TM},g,d,\tau)=\widehat{A}(TM,\nabla^{TM}){\rm ch}(\Theta_2(TM),\nabla^{\Theta_2(TM)},\tau){\rm ch}(Q_2(E),g^{Q_2(E)},d,\tau);\eqno(3.13)$$
$$\Phi_W'(\nabla^{TM},g,d,\tau)=\widehat{A}(TM,\nabla^{TM}){\rm ch}(\Theta_3(TM),\nabla^{\Theta_2(TM)},\tau) {\rm ch}(Q_3(E),g^{Q_3(E)},d,\tau).\eqno(3.14)$$
By Proposition 2.4 and Theorem 2.6 in [HY], we have that if $c_3(E_C,g,d)=0$, then for any integer $l\geq 1$ and $j=1,2,3$,~
$\Phi_L(\nabla^{TM},g,d,\tau)^{(4l-1)},~\Phi_W(\nabla^{TM},g,d,\tau)^{(4l-1)}$ and $\Phi_W'(\nabla^{TM},g,d,\tau)^{(4l-1)}$
are modular forms of weight $2l$ over $\Gamma_0(2)$, $\Gamma^0(2)$ and $\Gamma_\theta$ respectively. We have
$$\left<\Phi_L(\nabla^{TM},g,d,\tau),[M]\right>=-{\rm Ind}(T\otimes \triangle(TM)\otimes \Theta_1(TM)\otimes(Q_1(E),g^{Q_1(E)}));$$
$$\left<\Phi_W(\nabla^{TM},g,d,\tau),[M]\right>=-{\rm Ind}(T\otimes \Theta_2(TM)\otimes(Q_2(E),g^{Q_2(E)}));$$
$$\left<\Phi_W'(\nabla^{TM},g,d,\tau),[M]\right>=-{\rm Ind}(T\otimes \Theta_3(TM)\otimes(Q_3(E),g^{Q_3(E)})),\eqno(3.15)$$
where ${\rm Ind}(T\otimes \cdots)$ denotes the index of the Toeplitz operator. Let $\{\pm 2\pi\sqrt{-1}x_j\}$ for $1\leq j\leq 2r-1$ be the Chern roots of $TM\otimes {\bf C}$. Similar to the computations in [Li1], we have
$$\Phi_L(\nabla^{TM},g,d,\tau)=
2^{2r-1}\left(\prod_{j=1}^{2r-1}\frac{x_j\theta'(0,\tau)}{\theta(x_j,\tau)}
\frac{\theta_1(x_j,\tau)}{\theta_1(0,\tau)}\right) {\rm ch}(Q_1(E),g^{Q_1(E)},d,\tau);\eqno(3.16)$$
$$\Phi_W(\nabla^{TM},g,d,\tau)=
\left(\prod_{j=1}^{2r-1}\frac{x_j\theta'(0,\tau)}{\theta(x_j,\tau)}
\frac{\theta_2(x_j,\tau)}{\theta_2(0,\tau)}\right) {\rm ch}(Q_2(E),g^{Q_2(E)},d,\tau);\eqno(3.17)$$
$$\Phi_W'(\nabla^{TM},g,d,\tau)=
\left(\prod_{j=1}^{2r-1}\frac{x_j\theta'(0,\tau)}{\theta(x_j,\tau)}
\frac{\theta_3(x_j,\tau)}{\theta_3(0,\tau)}\right) {\rm ch}(Q_3(E),g^{Q_3(E)},d,\tau).\eqno(3.18)$$
Clearly,
$\Theta_1(T_{C}M)\otimes Q_1(E)$ and $\Theta_2(T_{C}M)\otimes Q_2(E)$ admit
formal Fourier expansion in $q^{\frac{1}{2}}$ as
$$\Theta_1(T_{C}M)\otimes Q_1(E)=A_0(T_{C}M,E)+A_1(T_{C}M,E)q
+\cdots,$$
$$\Theta_2(T_{C}M)\otimes Q_2(E)=B_0(T_{C}M,E)+B_1(T_{C}M,E)q
^{\frac{1}{2}}+\cdots,\eqno(3.19)$$ where the $A_j$ and $B_j$ are
elements in the semi-group formally generated by Hermitian vector
bundles over $M$. Moreover, they carry canonically induced Hermitian
connections. If $$B_j(T_{C}M,E)=B_{j,1}(T_{C}M)\otimes B_{j,2}(E),$$ we let
$$\widetilde{{\rm
ch}}(B_j(T_{C}M,E))={\rm ch}( B_{j,1}(T_{C}M)){\rm ch}(B_{j,2}(E),g^{B_{j,2}(E)},d).$$
If $\omega$ is a differential form over $M$, we denote
$\omega^{(4r-1)}$ its top degree component. Our main results include the following theorem.\\

\noindent {\bf Theorem 3.1} {\it If $c_3(E,g,d)=0$, then}
$$\left\{{\widehat{L}(TM,\nabla^{TM})}{\rm ch}(\triangle(E),g^{\triangle(E)},d)\right\}^{(4r-1)}
=2^{3r-1+\frac{N}{2}}\sum_{l=1}^{[\frac{r}{2}]}2^{-6l}h_l ,\eqno(3.20)$$ {\it
where each $h_l$,~$1\leq l\leq [\frac{r}{2}]$, is a canonical
integral linear combination of
$$\left\{\widehat{A}(TM,\nabla^{TM})\widetilde{{\rm
ch}}(B_j(T_{C}M,E))\right\}^{(4r-1)},$$  $1\leq j\leq
l$ and $h_1,h_2$ are given by (3.25) and (3.30).}
\\

\noindent {\bf Proof.} Similarly to the computations in [Li1, P.35] and by (2.26) in [HY] and the condition $c_3(E,g,d)=0$, we have
$$\Phi_W(\nabla^{TM},g,d,-\frac{1}{\tau})^{(4r-1)}=\frac{\tau^{2r}}{2^{2r-1+\frac{N}{2}}}\Phi_L(\nabla^{TM},g,d,\tau)^{(4r-1)}.\eqno(3.21)$$
By $\Phi_W(\nabla^{TM},g,d,\tau)^{(4r-1)}$ is a modular form of weight $2r$ over $\Gamma^0(2)$. By Lemma 2.2, we have
$$\Phi_W(\nabla^{TM},g,d,{\tau})^{(4r-1)}=h_0(8\delta_2)^{r}+h_1(8\delta_2)^{r-2}\varepsilon_2+\cdots+h_{[\frac{r}{2}]}(8\delta_2)
^{r-2[\frac{r}{2}]}\varepsilon^{[\frac{r}{2}]}_2,\eqno(3.22)$$ where
each $h_l,~ 0\leq l\leq [\frac{r}{2}],$ is a real multiple of the
volume form at $x$. By (2.20) (3.21) and (3.22), we get
$$\Phi_L(\nabla^{TM},g,d,{\tau})^{(4r-1)}=2^{2r-1+\frac{N}{2}}\left[h_0(8\delta_1)^{r}+h_1(8\delta_1)^{r-2}\varepsilon_1+\cdots+h_{[\frac{r}{2}]}(8\delta_1)
^{r-2[\frac{r}{2}]}\varepsilon_1^{[\frac{r}{2}]}\right].\eqno(3.23)$$
By comparing the constant term in (3.23), we get (3.20). By comparing
the coefficients of $q^{\frac{j}{2}}$,~$j\geq 0$ between the two
sides of (3.22), we can use the induction method to prove that $h_0=0$ and each
$h_l~1\leq l\leq [\frac{r}{2}]$, can be expressed through a
canonical integral linear combination of $$\left\{\widehat{A}(TM,\nabla^{TM})\widetilde{{\rm
ch}}(B_j(T_{C}M,E))\right\}^{(4r-1)},$$  $1\leq j\leq
l$. Direct computations shows that
$$\Theta_2(T_{C}M)\otimes Q_2(E)=1-(\widetilde{T_CM}+\widetilde{E_C})q^{\frac{1}{2}}+(\widetilde{T_CM}+\wedge^2\widetilde{T_CM}+\wedge^2\widetilde{E_C}+
\widetilde{T_CM}\otimes\widetilde{E_C})q+\cdots.\eqno(3.24)$$
By (2.19) and comparing the coefficient of $q^{\frac{1}{2}}$ of (3.22), we get
$$h_1=(-1)^{r-2}\left\{\widehat{A}(TM,\nabla^{TM})\widetilde{{\rm
ch}}(B_1(T_{C}M,E))\right\}^{(4r-1)}$$
$$=(-1)^{r-1}\left\{\widehat{A}(TM,\nabla^{TM}){\rm
ch}(E,g,d)\right\}^{(4r-1)}.\eqno(3.25)$$
By (2.19) and comparing the coefficient of $q$ of (3.22), we get
$$h_2=(-1)^{r-4}\left\{\widehat{A}(TM,\nabla^{TM})\widetilde{{\rm
ch}}(B_2(T_{C}M,E))\right\}^{(4r-1)}+[-8+24(-1)^r(r-2)]h_1.\eqno(3.26)$$
By (2.3), then
$$\widetilde{{\rm
ch}}(B_2(T_{C}M,E))={\rm ch}(\wedge^2\widetilde{E_C},g,d)+{\rm ch}(\widetilde{T_CM}){\rm ch}(E,g,d);\eqno(3.27)$$
$$\wedge^2\widetilde{E_C}=S^2({\bf C}^N)+\wedge^2E_C-E_C\otimes{\bf C}^N;\eqno(3.28)$$
$${\rm ch}(\wedge^2\widetilde{E_C},g,d)={\rm ch}(\wedge^2{E_C},g,d)-N{\rm ch}(E_C,g,d).\eqno(3.29)$$
By (3.25)-(3.29), we have
$$h_2=(-1)^{r-4}\left[\widehat{A}(TM,\nabla^{TM}){\rm ch}(T_CM){\rm ch}(E,g,d)+\widehat{A}(TM,\nabla^{TM}){\rm ch}(\wedge^2{E},g,d)\right]^
{(4r-1)}$$
$$
+[(-1)^{r-4}(7-N-4r)-24(r-2)+8(-1)^r]\left[\widehat{A}(TM,\nabla^{TM}){\rm ch}(E,g,d)\right]^
{(4r-1)}.\eqno(3.30)$$
$\Box$\\

\noindent {\bf Corollary 3.2} {\it Let $M$ be a $(4r-1)$-dimensional spin manifold and $c_3(E,g,d)=0$, then}
$${\rm Ind}(T\otimes\triangle{TM}\otimes(\triangle(E),g^{\triangle(E)}))
=-2^{3r-1+\frac{N}{2}}\sum_{l=1}^{[\frac{r}{2}]}2^{-6l}h_l ,\eqno(3.31)$$ {\it
where each $h_l$,~$1\leq l\leq [\frac{r}{2}]$, is a canonical
integral linear combination of
${\rm Ind}(T\otimes(B_j(T_{C}M,E)))$.}\\

\noindent {\bf Corollary 3.3} {\it Let $M$ be a $(4r-1)$-dimensional spin manifold and $c_3(E,g,d)=0$. If $r$ is even, then
$${\rm Ind}(T\otimes\triangle{TM}\otimes(\triangle(E),g^{\triangle(E)}))\equiv 0~~({\rm mod}~ 2^{\frac{N}{2}-1}).$$ If $r$ is odd, then
$${\rm Ind}(T\otimes\triangle{TM}\otimes(\triangle(E),g^{\triangle(E)}))\equiv 0~~({\rm mod}~ 2^{\frac{N}{2}+2}).$$}\\

\noindent {\bf Corollary 3.4} {\it Let $M$ be a $11$-dimensional manifold and $c_3(E,g,d)=0$. Then we have}
$$\left\{{\widehat{L}(TM,\nabla^{TM})}{\rm ch}(\triangle(E),g^{\triangle(E)},d)\right\}^{(11)}=2^{\frac{N}{2}+2}
\left\{{\widehat{A}(TM,\nabla^{TM})}{\rm ch}(E,g,d)\right\}^{(11)}.\eqno(3.32)$$\\

\noindent {\bf Corollary 3.5} {\it Let $M$ be a $15$-dimensional manifold and $c_3(E,g,d)=0$. Then we have}
$$\left\{\widehat{L}(TM,\nabla^{TM}){\rm ch}(\triangle(E),g^{\triangle(E)},d)\right\}^{(15)}=
2^{\frac{N}{2}-1}\left[-(113+N)
{\widehat{A}(TM,\nabla^{TM})}{\rm ch}(E,g,d)\right.$$
$$\left.+
{\widehat{A}(TM,\nabla^{TM})}{\rm ch}(T_CM){\rm ch}(E,g,d)+{\widehat{A}(TM,\nabla^{TM})}{\rm ch}(\wedge^2E,g,d)\right]^{(15)}.\eqno(3.33)$$\\

\indent By (5.21) in [CH], we have
\begin{eqnarray}
\Theta_1(TM)&=&2+2(T_CM-4r+1)q
+2[(-8r+3)T_CM\nonumber\\
&&+T_CM\otimes T_CM+(4r-1)(4r-2)]q^2+\cdots.\nonumber~~~~~~~~~~~~~~~~~~~~~~~~~~~~(3.34)
\end{eqnarray}
$${\rm ch}(Q_1(E),g^{Q_1(E)},d)={\rm ch}(\triangle(E)\otimes (E_Cq+((2-4r)E_C+\wedge^2E_C)q^2+\cdots),g,d).\eqno(3.35)$$
\begin{eqnarray}
\Theta_1(TM)\otimes Q_1(E)&=&2\triangle(E)\otimes E_Cq
+\left[2\triangle(E)\otimes ((2-4r)E_C+\wedge^2E_C)\right.\nonumber\\
&&\left.+2(T_CM-4r+1)\otimes \triangle(E)\otimes E_C\right]q^2+O(q^3).~~~~~~~~~~~~\nonumber(3.36)
\end{eqnarray}
By (5.22) in [CH], we have for $1\leq l \leq [\frac{r}{2}]$
\begin{eqnarray}
(8\delta_1)^{r-2l}\varepsilon^l_1
&=&2^{r-6l}\left[1+(24r-64l)q+(288r^2-1536rl\right.\nonumber\\
&&\left.+2048l^2+512l-264r)q^2+O(q^3)\right].~~~~~~~~~~~~~~~~~~~~~~~~~~~~~~~~~\nonumber(3.37)
\end{eqnarray}
By comparing the coefficient of $q$ in (3.23), we get\\

\noindent {\bf Theorem 3.6} {\it If $c_3(E,g,d)=0$, then}
\begin{eqnarray}
&&\left\{{\widehat{L}(TM,\nabla^{TM})}{\rm ch}(\triangle(E)\otimes E_C,g,d)\right.\left.-12r{\widehat{L}(TM,\nabla^{TM})}{\rm ch}(\triangle(E),g^{\triangle(E)},d)\right\}^{(4r-1)}\nonumber\\
&=&-2^{3r+4+\frac{N}{2}}\sum_{l=1}^{[\frac{r}{2}]}2^{-6l}lh_l . \nonumber~~~~~~~~~~~~~~~~~~~~~~~~~~~~~~~~~~~~~~~~~~~~~~~~~~~~~~~~~~~~~~~~~~~~(3.38)
\end{eqnarray}\\

\noindent {\bf Corollary 3.7} {\it Let $M$ be a $(4r-1)$-dimensional spin manifold and $c_3(E,g,d)=0$.\\ If $r$ is even, then
$${\rm Ind}(T\otimes\triangle(TM)\otimes(\triangle(E)\otimes E_C,g))-12r{\rm Ind}(T\otimes\triangle(TM)\otimes(\triangle(E),g))\equiv 0~~({\rm mod}~ 2^{\frac{N}{2}+4}).\eqno(3.39)$$
 If $r$ is odd, then}
$${\rm Ind}(T\otimes\triangle(TM)\otimes(\triangle(E)\otimes E_C,g))-12r{\rm Ind}(T\otimes\triangle(TM)\otimes(\triangle(E),g))\equiv 0~~({\rm mod}~ 2^{\frac{N}{2}+7}).\eqno(3.40)$$\\

By comparing the coefficient of $q^2$ in (3.23), we get\\

\noindent {\bf Theorem 3.8} {\it If $c_3(E,g,d)=0$, then}
\begin{eqnarray}
&&\left\{{\widehat{L}(TM,\nabla^{TM})}\left[(19-56r){\rm ch}(\triangle(E)\otimes E_C,g,d)-144(r^2-r){\rm ch}(\triangle(E),g,d)\right.\right.\nonumber\\
&&\left.\left.
+{\rm ch}(T_CM){\rm ch}(\triangle(E)\otimes E_C,g,d)+{\rm ch}(\triangle(E)\otimes \wedge^2E_C,g,d)\right]\right\}^{(4r-1)}\nonumber\\
&&=2^{3r+9+\frac{N}{2}}\sum_{l=1}^{[\frac{r}{2}]}2^{-6l}l^2h_l.\nonumber~~~~~~~~~~~~~~~~~~~~~~~~~~~~~~~~~~~~~~~~~~~~~~~~~~~~~~~~~~~~~~~~~~~~(3.41)
\end{eqnarray}\\

\noindent {\bf Corollary 3.9} {\it Let $M$ be a $(4r-1)$-dimensional spin manifold and $c_3(E,g,d)=0$.
Write $A$ for the index of the Toeplitz operator determined by the left hand of (3.38). If $r$ is even, then
$A\equiv 0~~({\rm mod}~ 2^{\frac{N}{2}+9}).$
 If $r$ is odd, then
$A\equiv 0~~({\rm mod}~ 2^{\frac{N}{2}+12}).$}\\

\noindent {\bf Corollary 3.10} {\it If ${\rm dim}M=11$ and $c_3(E,g,d)=0$, then}
\begin{eqnarray}&&\left\{{\widehat{L}(TM,\nabla^{TM})}{\rm ch}(\triangle(E)\otimes E_C,g,d)-36{\widehat{L}(TM,\nabla^{TM})}{\rm ch}(\triangle(E),g,d)\right.\nonumber\\
&&\left.
+2^{7+\frac{N}{2}}{\widehat{A}(TM,\nabla^{TM})}{\rm ch}( E,g,d)\right\}^{(11)}=0 .\nonumber~~~~~~~~~~~~~~~~~~~~~~~~~~~~~~~~~~~~~~~~~~~~(3.42)
\end{eqnarray}\\

\noindent {\bf Corollary 3.11} {\it If ${\rm dim}M=11$ and $c_3(E,g,d)=0$, then}
\begin{eqnarray}&&\left\{{\widehat{L}(TM,\nabla^{TM})}\left[-149{\rm ch}(\triangle(E)\otimes E_C,g,d)-864{\rm ch}(\triangle(E),g,d)\right.\right.\nonumber\\
&&\left.\left.
+{\rm ch}(T_CM){\rm ch}(\triangle(E)\otimes E_C,g,d)+{\rm ch}(\triangle(E)\otimes \wedge^2E_C,g,d)\right]\right\}^{(11)}\nonumber\\
&&
=2^{12+\frac{N}{2}}[{\widehat{A}(TM,\nabla^{TM})}{\rm ch}( E,g,d)]^{(11)}.\nonumber~~~~~~~~~~~~~~~~~~~~~~~~~~~~~~~~~~~~~~~~~~~~~~~~~(3.43)
\end{eqnarray}\\

\section{Twisted cancellation formulas in odd dimension}

   \quad Let $M$ be a $4r-1$ dimensional Riemannian
manifold and
 $\xi$ be a rank two real oriented Euclidean vector
   bundle over $M$ carrying with a Euclidean connection
   $\nabla^\xi$. Set
   $$\Theta_1(T_{C}M,\xi_C)=\bigotimes _{n=1}^{\infty}S_{q^n}(\widetilde{T_CM})\otimes
\bigotimes
_{m=1}^{\infty}\wedge_{q^m}(\widetilde{T_CM}-2\widetilde{\xi_C})\otimes
\bigotimes _{r=1}^{\infty}\wedge
_{q^{r-\frac{1}{2}}}(\widetilde{\xi_C})\otimes\bigotimes
_{s=1}^{\infty}\wedge _{-q^{s-\frac{1}{2}}}(\widetilde{\xi_C}),$$
$$\Theta_2(T_{C}M,\xi_C)=\bigotimes _{n=1}^{\infty}S_{q^n}(\widetilde{T_CM})\otimes
\bigotimes
_{m=1}^{\infty}\wedge_{-q^{m-\frac{1}{2}}}(\widetilde{T_CM}-2\widetilde{\xi_C})\otimes
\bigotimes _{r=1}^{\infty}\wedge
_{q^{r-\frac{1}{2}}}(\widetilde{\xi_C})\otimes\bigotimes
_{s=1}^{\infty}\wedge _{q^{s}}(\widetilde{\xi_C}),$$
$$\Theta_3(T_{C}M,\xi_C)=\bigotimes _{n=1}^{\infty}S_{q^n}(\widetilde{T_CM})\otimes
\bigotimes
_{m=1}^{\infty}\wedge_{q^{m-\frac{1}{2}}}(\widetilde{T_CM}-2\widetilde{\xi_C})\otimes
\bigotimes _{r=1}^{\infty}\wedge
_{q^{r}}(\widetilde{\xi_C})\otimes\bigotimes _{s=1}^{\infty}\wedge
_{-q^{s-\frac{1}{2}}}(\widetilde{\xi_C}).\eqno(4.1)$$ Let
$c=e(\xi,\nabla^{\xi})$ be the Euler form of $\xi$ canonically
associated to $\nabla^\xi$. Set
$$\Phi_L(\nabla^{TM},\nabla^{\xi},g,d,\tau)=\frac{\widehat{L}(TM,\nabla^{TM})}{{\rm
cosh}^2(\frac{c}{2})}{\rm
ch}(\Theta_1(T_{C}M,\xi_C),\nabla^{\Theta_1(T_{C}M,\xi_C)})$$
$$\cdot{\rm ch}(Q_1(E),g^{Q_1(E)},d,\tau),$$
$$\Phi_W(\nabla^{TM},\nabla^{\xi},d,g,\tau)={\widehat{A}(TM,\nabla^{TM})}{\rm
cosh}(\frac{c}{2}){\rm
ch}(\Theta_2(T_{C}M,\xi_C),\nabla^{\Theta_2(T_{C}M,\xi_C)})$$
$$\cdot {\rm ch}(Q_2(E),g^{Q_1(E)},d,\tau),$$
$$\Phi_W'(\nabla^{TM},\nabla^{\xi},d,g,\tau)={\widehat{A}(TM,\nabla^{TM})}{\rm
cosh}(\frac{c}{2}){\rm
ch}(\Theta_3(T_{C}M,\xi_C),\nabla^{\Theta_3(T_{C}M,\xi_C)})$$
$$\cdot {\rm ch}(Q_3(E),g^{Q_1(E)},d,\tau).\eqno(4.2)$$
Let $\{\pm2\pi\sqrt{-1}x_j| ~1\leq j\leq 2r-1\}$ and
$\{\pm2\pi\sqrt{-1}u\}$ be the Chern roots of $T_CM$ and $\xi_C$
respectively and $c=2\pi\sqrt{-1}u.$ Through direct computations, we
get (cf. [HZ2])
$$\Phi_L(\nabla^{TM},\nabla^{\xi},g,d,\tau)={2}^{2r-1}\left\{
\left(\prod_{j=1}^{2r-1}x_j\frac{\theta'(0,\tau)}{\theta(x_j,\tau)}\frac
{\theta_1(x_j,\tau)}{\theta_1(0,\tau)}\right)\frac{\theta_1^2(0,\tau)}
{\theta_1^2(u,\tau)}\frac{\theta_3(u,\tau)}{\theta_3(0,\tau)}
\frac{\theta_2(u,\tau)}{\theta_2(0,\tau)}\right\}$$
$$\cdot{\rm ch}(Q_1(E),g^{Q_1(E)},d,\tau);\eqno(4.3)$$
$$\Phi_W(\nabla^{TM},\nabla^{\xi},\tau)=\left(
\prod_{j=1}^{2r-1}x_j\frac{\theta'(0,\tau)}{\theta(x_j,\tau)}\frac
{\theta_2(x_j,\tau)}{\theta_2(0,\tau)}\right)\frac{\theta_2^2(0,\tau)}
{\theta_2^2(u,\tau)}\frac{\theta_3(u,\tau)}{\theta_3(0,\tau)}
\frac{\theta_1(u,\tau)}{\theta_1(0,\tau)}$$
$$\cdot{\rm ch}(Q_1(E),g^{Q_1(E)},d,\tau);\eqno(4.4)$$
$$\Phi_W'(\nabla^{TM},\nabla^{\xi},\tau)=\left(
\prod_{j=1}^{2r-1}x_j\frac{\theta'(0,\tau)}{\theta(x_j,\tau)}\frac
{\theta_3(x_j,\tau)}{\theta_3(0,\tau)}\right)\frac{\theta_3^2(0,\tau)}
{\theta_3^2(u,\tau)}\frac{\theta_1(u,\tau)}{\theta_1(0,\tau)}
\frac{\theta_2(u,\tau)}{\theta_2(0,\tau)}$$
$$\cdot{\rm ch}(Q_1(E),g^{Q_1(E)},d,\tau).\eqno(4.5)$$
Similarly to (3.21), we have
$$\Phi_W(\nabla^{TM},\nabla^{\xi},g,d,-\frac{1}{\tau})^{(4r-1)}
=\frac{\tau^{2r}}{2^{2r-1+\frac{N}{2}}}\Phi_L(\nabla^{TM},\nabla^{\xi},g,d,\tau)^{(4r-1)}.\eqno(4.6)$$
Similarly to Theorem 2.6 in [HY], we have\\

\noindent{\bf Proposition 4.1}~ {\it If $c_3(E_C,g,d)=0$, then for any integer $l\geq 1$ and $j=1,2,3$,~
$$\Phi_L(\nabla^{TM},\nabla^{\xi},g,d,\tau)^{(4l-1)},~\Phi_W(\nabla^{TM},\nabla^{\xi},g,d,\tau)^{(4l-1)}$$ and $\Phi_W'(\nabla^{TM},\nabla^{\xi},g,d,\tau)^{(4l-1)}$
are modular forms of weight $2l$ over $\Gamma_0(2)$, $\Gamma^0(2)$ and $\Gamma_\theta$ respectively.}\\

We know that (3.22) and (3.23) hold in the twisted case. We define $B_j(T_{C}M,\xi_C,E)$ similarly to $B_j(T_{C}M,E)$. Similarly to Theorem 3.1, we have\\

\noindent {\bf Theorem 4.2} {\it If $c_3(E,g,d)=0$, then}
$$\left\{\frac{\widehat{L}(TM,\nabla^{TM})}{{\rm
cosh}^2(\frac{c}{2})}{\rm ch}(\triangle(E),g^{\triangle(E)},d)\right\}^{(4r-1)}
=2^{3r-1+\frac{N}{2}}\sum_{l=1}^{[\frac{r}{2}]}2^{-6l}\widetilde{h_l} ,\eqno(4.7)$$ {\it
where each $\widetilde{h_l}$,~$1\leq l\leq [\frac{r}{2}]$, is a canonical
integral linear combination of
$$\left\{\widehat{A}(TM,\nabla^{TM}){\rm cosh}(\frac{c}{2})\widetilde{{\rm
ch}}(B_j(T_{C}M,\xi_C,E))\right\}^{(4r-1)},$$  $1\leq j\leq
l$.}\\

By the definition, we have
\begin{eqnarray}
&&\Theta_2(TM,\xi_C)\otimes Q_2(E)\nonumber\\&=&1+(3\widetilde{\xi_C}-\widetilde{T_CM}-\widetilde{E_C})q^{\frac{1}{2}}+\left[-3\widetilde{\xi_C}\otimes(\widetilde{T_CM}+\widetilde{E_C})\right.\nonumber\\
&&+(\widetilde{T_CM}+\wedge^2\widetilde{T_CM}+\wedge^2\widetilde{E_C}+
\widetilde{T_CM}\otimes\widetilde{E_C})\nonumber\\
&&\left.+(3\widetilde{\xi_C}\otimes \widetilde{\xi_C}+2S^2(\widetilde{\xi_C})+\wedge^2(\widetilde{\xi_C})+
\widetilde{\xi_C})\right]q+\cdots.\nonumber~~~~~~~~~~~~~~~~~~~~~~~~~~~~~~~~~~~~(4.8)
\end{eqnarray}
Comparing the coefficients of $q^{\frac{1}{2}}$ and $q$ in (3.22), we get
$$\widetilde{h_1}=(-1)^{r-1}\left\{\widehat{A}(TM,\nabla^{TM}){\rm cosh}\frac{c}{2}{\rm
ch}(E,g,d)\right\}^{(4r-1)}.\eqno(4.9)$$
\begin{eqnarray}
\widetilde{h_2}&=&(-1)^{r-4}\left[\widehat{A}(TM,\nabla^{TM}){\rm ch}(T_CM-3\xi_C){\rm ch}(E,g,d){\rm cosh}\frac{c}{2}\right.\nonumber\\
&&\left.+\widehat{A}(TM,\nabla^{TM}){\rm ch}(\wedge^2{E},g,d){\rm cosh}\frac{c}{2}\right]^
{(4r-1)}\nonumber\\
&&
+[(-1)^{r-4}(7-N-4r)-24(r-2)+8(-1)^r]\nonumber\\
&&\cdot\left[\widehat{A}(TM,\nabla^{TM}){\rm ch}(E,g,d){\rm cosh}\frac{c}{2}\right]^
{(4r-1)}.\nonumber~~~~~~~~~~~~~~~~~~~~~~~~~~~~~~~~~~~~~~~(4.10)
\end{eqnarray}\\

\noindent {\bf Corollary 4.3} {\it Let $M$ be a $(4r-1)$-dimensional spin$^c$ manifold and $c_3(E,g,d)=0$.
If $r$ is even, then
$$2^{1-\frac{N}{2}}\left<\frac{\widehat{L}(TM,\nabla^{TM})}{{\rm
cosh}^2(\frac{c}{2})}{\rm ch}(\triangle(E),g^{\triangle(E)},d),[M]\right>\equiv \widetilde{h}_{\frac{r}{2}}~~({\rm mod}~ 64{\bf Z}).\eqno(4.11)$$
 If $r$ is odd, then }
$$2^{-2-\frac{N}{2}}\left<\frac{\widehat{L}(TM,\nabla^{TM})}{{\rm
cosh}^2(\frac{c}{2})}{\rm ch}(\triangle(E),g^{\triangle(E)},d),[M]\right>\equiv \widetilde{h}_{\frac{r-1}{2}}~~({\rm mod}~ 64{\bf Z}).\eqno(4.12)$$\\

\noindent {\bf Corollary 4.4} {\it Let $M$ be a $11$-dimensional manifold and $c_3(E,g,d)=0$. Then we have}
$$\left\{\frac{\widehat{L}(TM,\nabla^{TM})}{{\rm
cosh}^2(\frac{c}{2})}{\rm ch}(\triangle(E),g^{\triangle(E)},d)\right\}^{(11)}=2^{\frac{N}{2}+2}
\left\{{\widehat{A}(TM,\nabla^{TM})}{\rm ch}(E,g,d){\rm cosh}\frac{c}{2}\right\}^{(11)}.\eqno(4.13)$$\\

\noindent {\bf Corollary 4.5} {\it Let $M$ be a $15$-dimensional manifold and $c_3(E,g,d)=0$. Then we have}
\begin{eqnarray}
&&\left\{\frac{\widehat{L}(TM,\nabla^{TM})}{{\rm
cosh}^2(\frac{c}{2})}{\rm ch}(\triangle(E),g^{\triangle(E)},d)\right\}^{(15)}\nonumber\\
&=&
2^{\frac{N}{2}-1}\left[-(113+N)
{\widehat{A}(TM,\nabla^{TM})}{\rm cosh}\frac{c}{2}{\rm ch}(E,g,d)\right.\nonumber\\
&&\left.+
{\widehat{A}(TM,\nabla^{TM})}{\rm ch}(T_CM-3\xi_C){\rm ch}(E,g,d){\rm cosh}\frac{c}{2}\right.\nonumber\\
&&\left.+{\widehat{A}(TM,\nabla^{TM})}{\rm ch}(\wedge^2E,g,d){\rm cosh}\frac{c}{2}\right]^{(15)}.\nonumber~~~~~~~~~~~~~~~~~~~~~~~~~~~~~~~~~~~~~~~~~~~~(4.14)
\end{eqnarray}\\

By the definition, we have
$$\widetilde{{\rm{ch}}}(\Theta_1(TM)\otimes Q_1(E))=2{\rm ch}(\triangle(E)\otimes E_C,g,d)q+O(q^{\frac{3}{2}}).\eqno(4.15)$$
Comparing the coefficient of $q$ in (3.23), we get\\

\noindent {\bf Theorem 4.6} {\it If $c_3(E,g,d)=0$, then}
\begin{eqnarray}
&&\left\{\frac{\widehat{L}(TM,\nabla^{TM})}{{\rm
cosh}^2(\frac{c}{2})}{\rm ch}(\triangle(E)\otimes E_C,g,d)\right.\left.-12r\frac{\widehat{L}(TM,\nabla^{TM})}{{\rm
cosh}^2(\frac{c}{2})}{\rm ch}(\triangle(E),g^{\triangle(E)},d)\right\}^{(4r-1)}\nonumber\\
&&=-2^{3r+4+\frac{N}{2}}\sum_{l=1}^{[\frac{r}{2}]}2^{-6l}l\widetilde{h_l}.\nonumber~~~~~~~~~~~~~~~~~~~~~~~~~~~~~~~~~~~~~~~~~~~~~~~~~~~~~~~~~~~~~~~~~~(4.16)
\end{eqnarray}\\

\noindent {\bf Corollary 4.7} {\it Let $M$ be a $(4r-1)$-dimensional spin$^c$ manifold and $c_3(E,g,d)=0$.
If $r$ is even, then
\begin{eqnarray}&&2^{-4-\frac{N}{2}}\left<\frac{\widehat{L}(TM,\nabla^{TM})}{{\rm
cosh}^2(\frac{c}{2})}{\rm ch}(\triangle(E)\otimes E_C,g,d)\right.\nonumber\\
&&\left.-12r\frac{\widehat{L}(TM,\nabla^{TM})}{{\rm
cosh}^2(\frac{c}{2})}{\rm ch}(\triangle(E),g^{\triangle(E)},d),[M]\right>\equiv -\frac{r}{2}\widetilde{h}_{\frac{r}{2}}~~({\rm mod}~ 64{\bf Z}).\nonumber~~~~~~~~~(4.17)
\end{eqnarray}
 If $r$ is odd, then }
\begin{eqnarray}&&2^{-6-\frac{N}{2}}\left<\frac{\widehat{L}(TM,\nabla^{TM})}{{\rm
cosh}^2(\frac{c}{2})}{\rm ch}(\triangle(E)\otimes E_C,g,d)\right.\nonumber\\
&&\left.-12r\frac{\widehat{L}(TM,\nabla^{TM})}{{\rm
cosh}^2(\frac{c}{2})}{\rm ch}(\triangle(E),g^{\triangle(E)},d),[M]\right>\equiv -\frac{r-1}{2}\widetilde{h}_{\frac{r-1}{2}}~~({\rm mod}~ 64{\bf Z}).\nonumber~~~(4.18)
\end{eqnarray}\\

\noindent {\bf Corollary 4.8} {\it~ If ${\rm dim}M=11$ and $c_3(E,g,d)=0$, then}
\begin{eqnarray}&&\left\{\frac{\widehat{L}(TM,\nabla^{TM})}{{\rm
cosh}^2(\frac{c}{2})}{\rm ch}(\triangle(E)\otimes E_C,g,d)-36\frac{\widehat{L}(TM,\nabla^{TM})}{{\rm
cosh}^2(\frac{c}{2})}{\rm ch}(\triangle(E),g,d)\right.\nonumber\\
&&\left.
+2^{7+\frac{N}{2}}{\widehat{A}(TM,\nabla^{TM})}{\rm ch}( E,g,d){\rm cosh}\frac{c}{2}\right\}^{(11)}=0 .\nonumber~~~~~~~~~~~~~~~~~~~~~~~~~~~~~~~~~~~~(4.19)
\end{eqnarray}\\

\noindent{\bf Remark.} 1. By the Dai-Zhang family index theorem associated to the Toeplitz operator in [DZ], similarly to [HL], we can easily
extend our anomaly cancellation formulas to the family case.\\
\indent 2. Similarly to [Li1], [HLZ] and [LW1], we also extend our anomaly cancellation formulas to these cases.\\

 \noindent{\bf Acknowledgement}~The
work of the first author was supported by NSF. The work of the second author was supported by NSFC. 11271062 and NCET-13-0721. The authors are indebted to Dr. F. Han for helpful comments. \\

\noindent {\bf References}\\

\noindent [AW] L. Alvarez-Gaum\'{e}, E. Witten, Graviational
anomalies,
Nucl. Phys. B234 (1983), 269-330.\\
 \noindent [Ch] K. Chandrasekharan, {\it Elliptic
Functions},
Spinger-Verlag, 1985. \\
\noindent [CH] Q. Chen, F. Han, Modular invariance and twisted
anomaly cancellations of characteristic numbers, Trans. Amer.
Math. Soc. 361 (2009), 1463-1493 \\
\noindent[CH1] Q. Chen, F. Han, Elliptic genera, transgression and loop space Chern-Simons forms. Comm. Anal. Geom. 17 (2009), no. 1, 73-106.\\
\noindent [DZ] X. Dai, W. Zhang, Higher spectral flow, J. Func. Anal., 157 (1998), 432-469.\\
\noindent[HL] F. Han, K. Liu, Gravitational Anomaly Cancellation and Modular Invariance, Alg. Geom. Topo. 14 (2014) 91-113.\\
\noindent[HLZ] F. Han, K. Liu and W. Zhang, Modular forms and generalized anomaly cancellation
formulas, J. Geom. Phys, 62 (2012) 1038-1053.\\
\noindent [HY] F. Han, J. Yu, On the Witten rigidity theorem for
odd dimensional manifolds, arXiv:1504.03007.\\
\noindent [HZ1] F. Han, W. Zhang, ${\rm
Spin}^c$-manifold and elliptic genera,
C. R. Acad. Sci. Paris Serie I., 336 (2003), 1011-1014.\\
\noindent [HZ2] F. Han, W. Zhang, Modular invariance, characteristic
numbers and eta Invariants, J.
Diff. Geom. 67 (2004), 257-288.\\
\noindent [Li1] K. Liu, Modular invariance and characteristic
numbers. Commu. Math. Phys. 174 (1995), 29-42.\\
\noindent[Li2] K. Liu, On elliptic genera and theta-functions.
Topology 35 (1996), no. 3, 617--640.\\
\noindent[Li3] K. Liu. On modular invariance and rigidity theorems. J. Differential
Geom., 41 (1995):343-396.\\
\noindent[LM] K. Liu and X. Ma. On family rigidity theorems. I. Duke
Math. J., 102 (2000):451-474.\\
\noindent[LMZ1] K. Liu, X. Ma, and W. Zhang. Spin$^c$ manifolds and
rigidity theorems in K-theory. Asian J. Math., 4 (2000):933-959.\\
\noindent[LMZ2] K. Liu, X. Ma, and W. Zhang. Rigidity and vanishing
theorems in K-theory. Comm. Anal. Geom., 11 (2003):121-180.\\
\noindent[LW] K. Liu and Y. Wang. Rigidity theorems on odd dimensional
manifolds. Pure Appl. Math. Q., 5 (2009):1139-1159.\\
\noindent[LW1] K. Liu and Y. Wang, A note on modular forms and generalized anomaly cancellation
formulas, Sci. China, Math., 56 (2013), 55-65.\\
\noindent [Zh] W. Zhang, {\it Lectures on Chern-weil Theory
and Witten Deformations.} Nankai Tracks in Mathematics Vol. 4, World
Scientific, Singapore, 2001.\\

\indent{ Center of Mathematical Sciences, Zhejiang University
Hangzhou Zhejiang 310027, China and Department of Mathematics,
University of California at Los Angeles, Los Angeles CA 90095-1555,
USA\\} \indent  Email: {\it liu@ucla.edu.cn; liu@cms.zju.edu.cn}\\

 \indent{School of Mathematics and Statistics,
Northeast Normal University, Changchun Jilin, 130024, China }\\
\indent E-mail: {\it wangy581@nenu.edu.cn }\\
\end {document}